\documentclass[leqno,12pt]{article}
\usepackage{latexsym}
\newtheorem{thm}{Theorem}[section]
\newtheorem{cor}[thm]{Corollary}

\newtheorem{note}[thm]{Note}

\newcommand{\ip}[2]{\mbox{$\langle #1,#2 \rangle$}}
\newcommand{\pf}{\noindent{\bf Proof\ \ }}
\newcommand{\finpf}{\hfill{$\Box$}\linespace}
\newcommand{\linespace}{\vspace{\baselineskip} \noindent}
\newcommand{\R}{{\bf R}}

\newcommand{\rt}{\rightarrow}
\def\tto{\;{\lower 1pt \hbox{$\rightarrow$}}\kern -12pt
           \hbox{\raise 2.8pt \hbox{$\rightarrow$}}\;}
\newenvironment{myequation}{\setcounter{equation}{\value{thm}}
   \begin{equation}}{\addtocounter{thm}{1}\end{equation}}

\newcommand{\bmye}{\begin{myequation}}
\newcommand{\emye}{\end{myequation}}

\begin{document}
\title{\textbf{\normalsize
THE STRUCTURED DISTANCE TO ILL-POSEDNESS FOR CONIC SYSTEMS}}
\author{\textbf{A.S. Lewis\thanks{Department 
of Mathematics, Simon Fraser University, Burnaby, BC V5A 1S6, Canada. 
\texttt{aslewis\char64 sfu.ca, 
http://www.cecm.sfu.ca/\~~\hspace{-4pt}aslewis}. 
Research supported by NSERC.} 
}}

\maketitle

\noindent {\bf Key words:}  condition number, conic system, distance
to infeasibility, structured singular values, sublinear maps,
surjectivity \\
{\bf AMS 2000 Subject Classification:} \\
Primary:  15A12, 90C31\\
Secondary: 65F35, 93B35

\begin{abstract}
An important measure of conditioning of a conic linear system is
the size of the smallest structured perturbation making the system
ill-posed.  We show that this measure is unchanged if we restrict to
perturbations of low rank.  We thereby derive a broad generalization of
the classical Eckart-Young result characterizing the distance to
ill-posedness for a linear map.
\end{abstract}

\section{Introduction}
Consider two finite-dimensional
normed spaces $X$ and $Y$, a fixed convex cone $K \subset X$, and
a linear mapping $A:X \rt Y$.  We call $A$ {\em well-posed}
if $AK=Y$.  In particular, in the purely linear
case $K=X$, well-posedness coincides with 
surjectivity.  Our interest is in the ``distance to ill-posedness'':
that is, we seek the smallest structured linear perturbation
$\Delta A : X \rt Y$ such that the perturbed mapping $A+\Delta A$ is 
not well-posed.  When $K=X$ and
the structure of perturbations is unrestricted, the classical
Eckart-Young theorem identifies the distance to ill-posedness as the
smallest singular value of $A$.

For more general convex cones $K$, and unstructured perturbations,
seminal work of Renegar \cite{Ren95a,Ren95} relates the distance to 
ill-posedness to the complexity of solving associated linear programs.
Imposing structure on the allowable perturbations (in order, for example,
to maintain a sparsity pattern in the map $A$) leads to a considerably more
involved theory.  In the purely linear case $K=X$, such questions arise as
``structured singular value'' calculations in the area of control theory,
pioneered by Doyle, known as ``$\mu$-analysis'' \cite{Doy82}.

In this article we follow quite closely the approach of Pe\~{n}a 
\cite{Pen03} in considering
structured perturbations to general conic systems.  We depend heavily
on the same rank-one reduction technique used in \cite{Pen03}
and introduced in \cite{Pen98,Pen00}.  Our approach differs
in several respects.  First, we develop
the theory in the concise and elegant language of sublinear set-valued
mappings (in other words, mappings whose graphs are convex cones).  This
notion substantially generalizes the idea of a conic convex system:
well-posedness becomes the notion of surjectivity of the mapping.
(In this framework, the unstructured case was developed in 
\cite{Lew99Ill}, and generalized in \cite{Don03}.)
Secondly, the structured perturbations we consider are rather general,
being of the form $\sum_i P_i T_i Q_i$ for linear mappings $T_i$
(where the linear mappings $P_i$ and $Q_i$ are fixed at the outset).  Thirdly,
we allow arbitrary norms on the underlying spaces.  Lastly, our proofs
consist of direct duality arguments, avoiding the necessity of ``lifting''
problems into higher dimensional spaces.  In this manner we hope to 
illuminate the structural simplicity of the key results.

The main result is as follows.  We consider finite-dimensional normed 
spaces $X,Y,U_i,V_i$, linear mappings $P_i:V_i \rt Y$ and 
$Q_i:X \rt U_i$ (for $i=1,2,\ldots,k$), and a surjective set-valued
mapping $F:X \tto Y$ with graph a closed convex cone.  Then, denoting
dual spaces and adjoint mappings by $*$, the
following four quantities are equal:
\[
\min_{\mbox{\scriptsize linear}~T_i}
\Big\{
\max_i \|T_i\| : F + \sum_i P_i T_i Q_i~\mbox{nonsurjective} 
\Big\};
\]
\[
\min_{\mbox{\scriptsize rank-one linear}~T_i}
\Big\{
\max_i \|T_i\| : F + \sum_i P_i T_i Q_i~\mbox{nonsurjective}
\Big\};
\]
\[
\min_{u_i^* \in U_i^*,\:z_i \ge 0,\:0 \ne y^* \in Y^*}
\Big\{
\max_i \frac{z_i}{\|P_i^* y_i\|} : 
\sum_i z_i Q_i^* u_i^* \in F^*(y^*),~\|u_i^*\| \le 1
\Big\};
\]
\[
\min_{v_i \in V_i,\:\|v_i\| \le 1}~
\sup_{x \in X,\:w_i > 0}
\Big\{
\min_i \frac{w_i}{\|Q_i x\|} : \sum_i w_i P_i v_i \in F(x)
\Big\}.
\]

\section{Rank-one perturbation}
As observed by Pe\~{n}a \cite{Pen98,Pen00}, the idea of rank-one pertubation is
fundamental to the theory of the distance to ill-posedness.  Our first, 
elementary result tries to capture the underlying idea in a way that extends
to structured perturbations.

Throughout this article we follow the terminology of \cite{Roc98}.  We call
a set-valued mapping $F: X \tto Y$ {\em positively-homogeneous} if its
{\em graph} 
\[
\mbox{gph}\,F = \{(x,y) \in X \times Y : y \in F(x)\}
\]
is a {\em cone} (which is to say, nonempty and closed under nonnegative
scalar multiplication).  To recapture the theory of conic linear 
systems we typically consider examples of the form
\[
F(x) = \left\{
\begin{array}{cl}
\{Ax\} & (x \in K)  \\
\emptyset & (x \not\in K),
\end{array}
\right.
\]
where the mapping $A:X \rt Y$ is linear and $K \subset X$ is a convex cone.
The {\em inverse} of a set-valued mapping $F$ is the mapping
$F^{-1}: Y \tto X$ defined by 
\[
x \in F^{-1}(y)~~\Leftrightarrow~~y \in F(x).
\]
We call $F$ {\em singular} if $F^{-1}(0) \ne \{0\}$.

We typically denote the norm on a normed space $X$ by $\| \cdot \|$
(or by $\|\cdot\|_X$ if we wish to be specific) and
the closed unit ball in $X$ by $B_X$, and we denote
the space of linear mappings from $X$ to $Y$ by $L(X,Y)$.
In particular, for a mapping $A \in L(X,Y)$, we denote the usual operator norm
by $\|A\|$.  We denote the dual space of $X$ by $X^*$, and we write
the action of a linear functional $x^* \in X^*$ on an
element $x \in X$ as $\ip{x^*}{x}$.  We are particularly interested
in {\em rank-one} mappings in $L(X,Y)$, which are those mappings of the
form $x \in X \mapsto \ip{x^*}{x}y$ for some given elements $x^* \in X^*$
and $y \in Y$:  we denote the set of such mappings by $L_1(X,Y)$.   
The norm of this mapping is just $\|x^*\| \cdot \|y\|$.

In what follows, we interpret $1/0 = +\infty$ and $1/{+\infty} = 0$.

\begin{thm}[rank-one reduction]  \label{2.1}
Consider finite-dimensional normed spaces $X,Y,U,V$, a positively-homogeneous set-valued
mapping $F:X \tto Y$, and linear mappings $P:V \rt Y$ and $Q:X \rt U$.
Then the quantity in $[0,+\infty]$ defined by
\[
\alpha = \inf_{T \in L(U,V)} \Big\{ \|T\| : F+PTQ~\mbox{singular} \Big\}
\]
is unchanged if we further restrict the infimum to be over mappings $T$ of 
rank one.  Furthermore, if we assume
\[
0 \in F(x)~\mbox{and}~x \ne 0~~\Rightarrow~~Qx \ne 0
\]
(as holds in particular if $Q$ is injective or $F$ is nonsingular),
then 
\[
\frac{1}{\alpha} = \sup_{x \in X,\:v \in B_V} \Big\{ \|Qx\| : Pv \in F(x) \Big\}.
\]
\end{thm}

\noindent
{\bf Note}~
We address the question of the attainment in the above infimum and supremum
in the next section.

\medskip

\pf
Denote the right hand side of the last equation by $\beta$.
Consider first the case where $F$ is singular.  In this case, clearly $\alpha = 0$,
and is attained by choosing the rank-one mapping $T=0$.  Choose any nonzero 
$x_1 \in F^{-1}(0)$, so by assumption, $Qx_1 \ne 0$.  Now by choosing
$x = \lambda x_1$ with $\lambda \in \R_+$ and $v=0$ in the definition
of $\beta$, and letting $\lambda$ grow,
we see $\beta = +\infty$, so the result holds.  We can therefore assume 
$F$ is nonsingular.

We next show $\alpha \ge 1/\beta$.  Consider any feasible mapping $T$ in the definition
of $\alpha$, so there exists a nonzero vector $x \in (F+PTQ)^{-1}(0)$.  Hence
we have $-PTQx \in F(x)$,  so since $F^{-1}(0) = \{0\}$, we deduce $TQx \ne 0$.
Positive homogeneity now implies 
\[
P\Big( -\frac{1}{\|TQx\|} TQx \Big) \in F\Big( \frac{1}{\|TQx\|} x \Big),
\]
so by definition, 
\[
\beta \ge \Big\| Q\|TQx\|^{-1}x \Big\| \ge \frac{1}{\|T\|}.
\]
Thus all feasible $T$ satisfy $\|T\| \ge \beta$, and we deduce $\alpha \ge 1/\beta$.

Next we define the quantity
\[
\gamma = \inf_{T \in L_1(U,V)} \Big\{ \|T\| : (F+PTQ)^{-1}(0) \neq \{0\} \Big\}.
\]
Clearly we have the inequality $\gamma \ge \alpha$, so it now suffices to
prove $\gamma \le 1/\beta$.  If $\beta = 0$ there is nothing to prove, so
we can assume $\beta > 0$.

Consider any feasible vectors $x$ and $v$ in the definition of $\beta$.
Since $\beta > 0$ we can assume $Qx \ne 0$.  There exists a norm-one
linear functional $u^* \in U^*$
satisfying $\ip{u^*}{Qx} = \|Qx\|$.  Now we have
\[
0 \in F(x) - Pv = F(x) - PTQx
\]
where $T:U \rt V$ is the rank-one linear map defined by
\[
Tu = \frac{\ip{u^*}{u}}{\|Qx\|} v.
\] 
Since we know $\|u\|_* = 1$ and $\|v\| \le 1$, we deduce
\[
\gamma \le \|T\| \le \frac{1}{\|Qx\|},
\]
so $1/\gamma \ge \|Qx\|$.  Finally, taking the supremum over all feasible 
vectors $x$ and $v$ in the definition of $\beta$ shows $1/\gamma \ge \beta$,
as required.
\finpf

\noindent
Notice that,  
if $X = Y$, the mapping $F$ is single-valued and linear,
and the mappings $P$ and $Q$ are just the identity, then we recover the classical
Eckart-Young theorem.

We next generalize to perturbations with a composite structure.  In conformity
with our previous usage, for $z \in \R_+$ we define
\[
\frac{z}{0} = 
\left\{
\begin{array}{cl}
+\infty & (z>0) \\
0       & (z=0).
\end{array}
\right.
\]

\begin{cor}[rank-one reduction for sums] \label{2.2}
Given finite-dimensional normed spaces $X,Y,U_i,V_i$, a 
positively-homogeneous set-valued 
mapping $F:X \tto Y$, and linear mappings $P_i:V_i \rt Y$ and $Q_i:X \rt U_i$
(for $i=1,2,\ldots,k$), the quantity
\[
\alpha := \inf_{T_i \in L(U_i,V_i)} \Big\{ \max_i \|T_i\| : 
F + \sum_i P_iT_iQ_i~\mbox{singular} \Big\}
\]
is unchanged if we further restrict the infimum to be over mappings $T_i$ of 
rank one.  Consequently we have the following:
\begin{eqnarray*}
\alpha 
& = &
\inf_{v_i \in B_{V_i},\:z_i \in \R_+,\:0 \ne x \in X}
\Big\{ \max_i \frac{z_i}{\|Q_i x\|} : \sum_i z_i P_i v_i \in F(x) \Big\}  \\
& = &
\inf_{v_i \in V_i,\: u_i^* \in B_{U_i^*},\: 0 \ne x \in X}
\Big\{
\max_i \|v_i\| : \sum_i \ip{u_i^*}{Q_i x} P_i v_i \in F(-x),  \\
& & \mbox{} \hspace{8cm} \ip{u_i^*}{Q_i x} \ge 0~\forall i
\Big\}.
\end{eqnarray*}
\end{cor}

\noindent
{\bf Note}~
As before, we address the question of the attainment in the above infima
in the next section.

\medskip

\pf
Fix any real $\epsilon > 0$ and consider any feasible mappings $T_i$ in the above
infimum.  By applying the preceding theorem we see there exists a
mapping $\hat{T}_k \in L_1(U_k,V_k)$ satisfying 
$\|\hat{T}_k\| < \|T_k\| + \epsilon$ and
\[
\Big( F + \sum_{i=1}^{k-1} P_iT_iQ_i + P_k \hat{T}_k Q_k \Big)^{-1}(0) \neq \{0\}.
\]
We can continue in this fashion, arriving at mappings
$\hat{T}_i \in L_1(U_i,V_i)$ satisfying 
$\|\hat{T}_i\| < \|T_i\| + \epsilon$ (for $i=1,2,\ldots,k$) and
\[
\Big( F + \sum_i P_i \hat{T}_i Q_i \Big)^{-1}(0) \neq \{0\}.
\]
Since $\epsilon > 0$ was arbitrary, the rank-one reduction now follows.

Consequently, we have $\alpha = \alpha_1$, where
\begin{eqnarray*}
\alpha_1 
&:=&
\inf_{v_i \in V_i,\:u_i^* \in U_i^*,\:0 \ne x \in X}
\Big\{ 
\max_i \|v_i\| \|u_i^*\| : \sum_i \ip{u_i^*}{Q_i x} P_i v_i \in F(x)
\Big\}  \\ 
&\le&
\inf_{v_i \in B_{V_i},\:u_i^* \in U_i^*,\:0 \ne x \in X}
\Big\{ 
\max_i \|v_i\| \|u_i^*\| : \sum_i \ip{u_i^*}{Q_i x} P_i v_i \in F(x)
\Big\}  \\
&\le&
\alpha_2,
\end{eqnarray*}
where
\[
\alpha_2 := 
\inf_{v_i \in B_{V_i},\:u_i^* \in U_i^*,\:0 \ne x \in X}
\Big\{ 
\max_i \|u_i^*\| : \sum_i \ip{u_i^*}{Q_i x} P_i v_i \in F(x)
\Big\}.
\]
On the other hand, suppose the vectors $v_i$, $u_i^*$ and $x$ are feasible
in the infimum defining $\alpha_1$.  If we define, for each index $i$,
\[
(\hat{v}_i , \hat{u}_i^*) = 
\left\{
\begin{array}{cl}
(\|v_i\|^{-1} v_i , \|v_i\| u_i^*) & (v_i \ne 0)  \\
(0,0)                            & (v_i = 0),
\end{array}
\right.
\]
then the vectors $\hat{v}_i$, $\hat{u}_i^*$ and $x$ are feasible in
the infimum defining $\alpha_2$, and $\|\hat{u}_i^*\| = \|v_i\| \|u_i^*\|$
for each $i$.  This proves $\alpha_2 \le \alpha_1$, so in fact 
$\alpha = \alpha_1 = \alpha_2$.

A completely analogous argument shows
\[
\alpha = 
\inf_{v_i \in V_i,\: u_i^* \in B_{U_i^*},\: 0 \ne x \in X}
\Big\{
\max_i \|v_i\| : \sum_i \ip{u_i^*}{Q_i x} P_i v_i \in F(-x)
\Big\}.
\]
The final expression for $\alpha$ claimed in the theorem now follows,
since the additional conditions $\ip{u_i^*}{Q_i x} \ge 0$ impose no essential
restriction:  for any index $i$ we can always replace the pair of vectors
$(v_i,u_i^*)$ with $(-v_i,-u_i^*)$ without changing feasibility or the
objective value.

Considering the definition of $\alpha_2$, we observe, for any vectors $v_i$,
\begin{eqnarray*}
\lefteqn{
\inf_{u_i^* \in U_i^*,\:0 \ne x \in X}
\Big\{ 
\max_i \|u_i^*\| : \sum_i \ip{u_i^*}{Q_i x} P_i v_i \in F(x)
\Big\}  }  \\
&=&
\inf_{u_i^* \in U_i^*,\:0 \ne x \in X,\:z_i \in \R_+}
\Big\{ 
\max_i \|u_i^*\| : \sum_i z_i P_i v_i \in F(x),~\ip{u_i^*}{Q_i x} = z_i,
\Big\}
\end{eqnarray*}
since a feasible choice of the variables on the right hand side immediately
gives a feasible choice on the left hand side with the same objective
value, while for any feasible choice 
of vectors $u_i^*$ and $x$ on the left hand side, setting 
$\hat{u}_i^* = (\mbox{sgn}\ip{u_i^*}{Q_i x}) u_i^*$ and 
$\hat{z_i} = |\ip{u_i^*}{Q_i x}|$ for each index $i$ gives a feasible choice
on the right hand side with the same objective value.

By observing that, for any vector $x \in X$ and
scalar $z_i \in \R_+$, we have
\[
\inf_{u_i^* \in U_i^*} \Big\{ \|u_i^*\| : \ip{u_i^*}{Q_i x} = z_i \Big\}
= \frac{z_i}{\|Q_i x\|},
\]
the result now follows.
\finpf

\noindent
{\bf Note}~~
It is not hard to see that the case $k=1$ gives back Theorem \ref{2.1}.

\section{Duality and surjectivity}
We return to our motivating example of the well-posedness of a linear 
mapping $A:X \rt Y$ relative to a convex cone $K \subset X$ (by which
we mean $AK=Y$).  If, as before, we define an associated set-valued
mapping $F:X \tto Y$ by
\bmye \label{3.1}
F(x) = \left\{
\begin{array}{cl}
\{Ax\} & (x \in K)  \\
\emptyset & (x \not\in K),
\end{array}
\right.
\emye
then well-posedness holds exactly when $F(X)=Y$.

We call a general set-valued mapping $F:X \tto Y$ {\em surjective} if
$F(X)=Y$, {\em closed} if its graph is closed, and {\em sublinear}
if its graph is a convex cone.  Sublinear set-valued mappings are
also known as {\em convex processes}.  The notions of singularity and 
surjectiveness are intimately connected via duality:  the {\em adjoint}
of $F$ is the set-valued mapping $F^* : Y^* \rt X^*$ defined by
\[
x^* \in F^*(y^*)~~\Leftrightarrow~~
\ip{y^*}{y} \ge \ip{x^*}{x}~\mbox{whenever}~y \in F(x).
\]
The adjoint is easily seen to be
closed and sublinear, and coincides with the
classical notion for single-valued linear mappings.  More generally,
direct calculation shows that for any linear mapping $G:X \rt Y$ we
have $(F+G)^* = F^* + G^*$.  It is simple to check that the adjoint 
of the set-valued mapping (\ref{3.1}) is defined by
$F^*(y^*) = A^* y^* + K^*$,
where $K^* \subset X^*$ is the usual (negative) polar cone for $K$.

The relationship between surjectiveness and singularity is described
by the following concise result, a special case of an infinite-dimensional
version of the open mapping theorem \cite{Bor86}.

\begin{thm}[open mapping] \label{3.2}
For finite-dimensional normed spaces $X$ and $Y$, a closed sublinear
set-valued mapping $F:X \tto Y$ is surjective if and only if its
adjoint mapping $F^*$ is nonsingular.
\end{thm}

\begin{note}
{\rm
If the closed sublinear set-valued mapping $F$ is surjective, then so is the
mapping $F+G$ for all small linear mappings $G$, and the analogous result
also holds for nonsingularity \cite{Rob76}.  Hence with
this assumption on $F$ in Theorem \ref{2.1} (rank-one reduction), the infimum 
\[
\inf_{T \in L(U,V)} \Big\{ \|T\| : F+PTQ~\mbox{singular} \Big\}
\]
is attained whenever finite, since it seeks the norm of the smallest element
in a nonempty closed set.  In this case, following the proof shows both the same
infimum over the rank-one mappings $T$ and the supremum
\[
\sup_{x \in X,\:v \in B_V} \Big\{ \|Qx\| : Pv \in F(x) \Big\}
\]
are also attained.
}
\end{note}

\begin{note} \label{3.4}
{\rm 
Using the preceding note, if the closed sublinear set-valued mapping $F$ is 
surjective in Corollary \ref{2.2} (rank-one reduction for sums), then the
infimum
\[
\inf_{T_i \in L(U_i,V_i)} \Big\{ \max_i \|T_i\| : 
F + \sum_i P_iT_iQ_i~\mbox{singular} \Big\}
\]
is attained whenever finite, whether over general or rank-one linear mappings $T_i$,
and in this case the infimum
\[
\inf_{v_i \in B_{V_i},\:z_i \in \R_+,\:0 \ne x \in X}
\Big\{ \max_i \frac{z_i}{\|Q_i x\|} : \sum_i z_i P_i v_i \in F(x) \Big\}
\]
is also attained.
}
\end{note}

Using the open mapping theorem (\ref{3.2}), we can quickly 
derive a version of Corollary \ref{2.2}
(rank-one reduction for sums) for nonsurjectivity rather than singularity.

\begin{thm}[rank reduction and surjectivity] \label{3.5n}
For any finite-dimen\-sional normed spaces $X,Y,U_i,V_i$, 
closed sublinear set-valued 
mapping $F:X \tto Y$, and linear mappings $P_i:V_i \rt Y$ and $Q_i:X \rt U_i$
(for $i=1,2,\ldots,k$), the quantity
\[
\alpha := \inf_{T_i \in L(U_i,V_i)} \Big\{ \max_i \|T_i\| : 
F + \sum_i P_iT_iQ_i~\mbox{nonsurjective} \Big\}
\]
is unchanged if we further restrict the infimum to be over mappings $T_i$ of 
rank one, and in fact
\begin{eqnarray*}
\alpha 
& = &
\inf_{u^*_i \in B_{U^*_i},\:z_i \in \R_+,\:0 \ne y^* \in Y^*}
\Big\{ \max_i \frac{z_i}{\|P^*_i y^*\|} : 
\sum_i z_i Q^*_i u^*_i \in F^*(y^*) \Big\} \\
& = &
\inf_{v_i \in B_{V_i},\: u_i^* \in U_i^*,\: 0 \ne y^* \in Y^*}
\Big\{
\max_i \|u_i^*\| : \sum_i \ip{y^*}{P_i v_i} Q_i^* u_i^* \in F^*(-y^*),  \\
& & \mbox{} \hspace{9cm} \ip{y^*}{P_i v_i} \ge 0~\forall i
\Big\}.
\end{eqnarray*}
Furthermore, all four infima are attained if $\alpha$ is finite.
\end{thm}

\pf
By the open mapping theorem, we have
\begin{eqnarray*}
\alpha
&=& 
\inf_{T_i \in L(U_i,V_i)} \Big\{ \max_i \|T_i\| : 
\Big( F + \sum_i P_iT_iQ_i \Big)^*~\mbox{singular} \Big\}  \\
&=& 
\inf_{T_i \in L(U_i,V_i)} \Big\{ \max_i \|T^*_i\| : 
F^* + \sum_i Q^*_iT^*_iP^*_i~\mbox{singular} \Big\},
\end{eqnarray*}
since the adjoint transformation $* : L(U_i,V_i) \rt L(V_i^*,U_i^*)$
leaves the norm fixed.  This transformation is in fact a bijection, which
also preserves the classes of rank-one mappings.  Corollary \ref{2.2} 
ensures the infimum is unchanged if we restrict to mappings $T_i$ for
which $T_i^*$ is rank-one, or in other words to rank-one $T_i$, as
required.  The final expressions follow directly from Corollary \ref{2.2}.
The final claim concerning attainment follows from Note \ref{3.4}.
\finpf

\section{Duality}
Our ultimate aim is to express the structured distance to nonsurjectivity 
in terms involving the mapping $F$ rather than its adjoint.  For this purpose,
the following result is crucial.

\begin{thm}[theorem of the alternative] \label{4.1}
For any finite-dimensional normed spaces $X,Y,U_i$, 
surjective closed sublinear set-valued 
mapping $F:X \tto Y$, linear mappings $Q_i:X \rt U_i$, and 
vectors $y_i \in Y$ (for $i=1,2,\ldots,k$),
exactly one of the following two systems has a solution:
\begin{enumerate}
\item[{\rm (}i\,{\rm )}]
\mbox{} \hfill
$\sum_i w_i y_i \in F(x)$,~~$\|Q_i x\| < w_i \in \R$ for each $i$,~~ 
$x \in X$;
\hfill \mbox{}
\item[{\rm (}ii\,{\rm )}] 
\mbox{} \hfill
$\sum_i \ip{y^*}{y_i}Q_i^* u_i^* \in F^*(-y^*)$,~~$0 \ne y^* \in Y^*$, 
\hfill \mbox{} \\
\mbox{} \hfill
$\ip{y^*}{y_i} \ge 0$~~and~~$u_i^* \in B_{U_i^*}$ for each $i$.
\hfill \mbox{}
\end{enumerate}
\end{thm}

\pf
Suppose first that both systems have solutions.  By the definition of the
adjoint, we deduce the inequality
\[
\Big\langle
-y^* , \sum_i w_i y_i
\Big\rangle 
\ge 
\Big\langle
\sum_i \ip{y^*}{y_i}Q_i^* u_i^* , x
\Big\rangle
\]
or equivalently
\[
0 \ge \sum_i \ip{y^*}{y_i} \Big( w_i + \ip{u_i^*}{Q_i x} \Big).
\]
Now each term in the sum on the right hand side is a product of two
factors, the first of which is nonnegative and the second of which is
strictly positive.  Hence this inequality can only hold if
$\ip{y^*}{y_i} = 0$ for each index $i$, and in this case we deduce
$0 \in F^*(-y^*)$.  But the mapping $F$ is surjective, so by the open
mapping theorem (\ref{3.2}) its adjoint $F^*$ is nonsingular, and this
is a contradiction.  Hence at most one of the two systems has a solution.

Suppose now that system $(i)$ has no solution.  Then the two convex subsets
of $X \times \R^k$
\[ 
\Big\{ (x,w) : \sum_i w_i y_i \in F(x) \Big\}
~~\mbox{and}~~
\Big\{ (x,w) : \|Q_i x\| < w_i~\mbox{for each}~i \Big\}
\]
are disjoint.  Both sets are clearly nonempty, so there exists a 
separating hyperplane:  there exists a nonzero vector
$(x^*,w^*) \in X^* \times R^k$ and a real $\mu$ such that
the two implications
\begin{eqnarray*}
\sum_i w_i y_i \in F(x) 
 & \Rightarrow &
\ip{x^*}{x} - \sum_i w_i^* w_i \ge \mu  \\
\|Q_i x\| < w_i~\mbox{for each}~i 
 & \Rightarrow &
\ip{x^*}{x} - \sum_i w_i^* w_i \le \mu.
\end{eqnarray*}

Considering the first implication, 
by the positive homogeneity of $F$, we deduce 
\bmye \label{4.2}
\sum_i w_i y_i \in F(x) ~~\Rightarrow~~
\ip{x^*}{x} - \sum_i w_i^* w_i \ge 0.
\emye
and $\mu \le 0$.  This, in conjunction with the 
second implication, shows 
\bmye \label{4.3}
w_i^* \ge 0~~\mbox{for each}~i,
\emye
and
\[
\ip{x^*}{x} \le \sum_i w_i^* \|Q_i x\|~~\mbox{for all}~x \in X.
\]
This inequality expresses the fact that the vector $x^*$ is a subgradient
at the origin for the convex function
\[
x \mapsto \sum_i w_i^* \|Q_i x\|,
\]
so by standard convex analysis we deduce
\bmye \label{4.4}
x^* \in \sum_i w_i^* Q_i^* B_{U_i^*}.
\emye

We now apply a rather standard duality argument to the implication (\ref{4.2}).
We define a function $f:Y \rt [-\infty,+\infty]$ by
\[
f(y) = 
\inf_{x \in X,\:w_i \in \R} 
\Big\{
\ip{x^*}{x} - \sum_i w_i^* w_i : y + \sum_i w_i y_i \in F(x)
\Big\}.
\]
Implication (\ref{4.2}) shows $f(0) = 0$, and a standard elementary argument
using the convexity of the graph of $F$ shows $f$ is convex.  Since the mapping
$F$ is surjective, the function $f$ never takes the value $+\infty$.  Consequently
(see \cite{Roc70}), $f$ has a subgradient $y^* \in Y^*$ at the origin, or in other
words,
\[
y+\sum_i w_i y_i \in F(x)
~~\Rightarrow~~
\ip{y^*}{y} \le \ip{x^*}{x} - \sum_i w_i^* w_i.
\]
Setting $x=0$ and $y = -\sum_i w_i y_i$ shows
\[
\sum_i w_i \Big( w_i^* - \ip{y^*}{y_i} \Big) \le 0~~\mbox{for all}~w \in \R^k,
\]
so 
\bmye \label{4.5}
w_i^* = \ip{y^*}{y_i}~~\mbox{for each}~i.  
\emye
Furthermore, setting each $w_i = 0$
shows
\[
y \in F(x)
~~\Rightarrow~~
\ip{y^*}{y} \le \ip{x^*}{x},
\]
or in other words,
\bmye \label{4.6}
-x^* \in F^*(-y^*).
\emye
Finally, putting together the relationships (\ref{4.3}), (\ref{4.4}), 
(\ref{4.5}), and (\ref{4.6}), shows we have constructed a solution to
system $(ii)$ in the theorem statement, as required.
\finpf

A helpful restatement of the above theorem is contained in the following
duality result.  Recall our convention $z/0 = +\infty$ for real $z>0$.

\begin{thm}[duality] \label{4.7}
Consider finite-dimensional normed spaces $X$, $Y$, $U_i$, 
a surjective closed sublinear set-valued 
mapping $F:X \tto Y$, linear mappings $Q_i:X \rt U_i$, and 
vectors $y_i \in Y$ (for $i=1,2,\ldots,k$).  Then the 
function $\Phi : Y^k \rt [0,+\infty]$ defined by
\bmye \label{3.5}
\Phi\Big((y_i)\Big) = 
\sup_{x \in X,\:0 < w_i \in \R}
\Big\{ \min_i \frac{w_i}{\|Q_i x\|} : \sum_i w_i y_i \in F(x) \Big\}.
\emye
is lower semicontinuous, and 
\begin{eqnarray*}
\lefteqn{
\Phi\Big((y_i)\Big) = 
}  \\
 & & 
\inf_{u_i^* \in U_i^*,\:0 \ne y^* \in Y^*}
\Big\{
\max_i \|u_i^*\| : \sum_i \ip{y^*}{y_i} Q_i^* u_i^* \in F^*(-y^*),~
\ip{y^*}{y_i} \ge 0~\forall i
\Big\}.
\end{eqnarray*}
Furthermore, the infimum on the right hand side is attained whenever finite.
\end{thm}

\pf
We first prove the lower semicontinuity.
For each index $i$ consider a sequence of vectors $y_i^r \rt y_i$ in the 
space $Y$, and consider a sequence of reals $s^r \rt s$ as $r \rt \infty$
satisfying $s^r \ge \Phi((y_i^r))$, or in other words 
\bmye \label{3.7}
\mbox{}~~~~~~x \in X,~0 < w_i \in \R~\mbox{and}~\sum_i w_i y_i^r \in F(x)
~~\Rightarrow~~
s^r \ge \min_i \frac{w_i}{\|Q_i x\|}.
\emye
Consider reals $w_i >0$ (for each $i$) satisfing 
$\sum_i w_i y_i \in F(\bar{x})$.  We want to show the inequality
\[
s \ge \min_i \frac{w_i}{\|Q_i \bar{x}\|}.
\]

To see this, we first note that, since $F$ is surjective, it is everywhere
{\em open}:  the image under $F$ of any open set is open.  In particular,
for any real $\delta > 0$, the set $F(\bar{x} + \mbox{int}\,\delta B_X)$ is an
open neighbourhood of the vector $\sum_i w_i y_i$, so for large $r$
must contain the point $\sum_i w_i y_i^r$.  Using this tool, we see
there exists a subsequence $R$ of the natural numbers such that
\begin{eqnarray*}
\sum_i w_i y_i^r              & \in & F(x^r)~\mbox{for all}~r \in R  \\
\lim_{r \rt \infty,\:r \in R} x^r &  =  & \bar{x}.
\end{eqnarray*}
Applying property (\ref{3.7}) shows
\[
s^r \ge \min_i \frac{w_i}{\|Q_i x^r\|}~~\mbox{for all}~r \in R.
\]
Hence there exists an index
$j \in \{1,2,\ldots,k\}$ and a further subsequence $R'$ of $R$ such that
\[
s^r \ge \frac{w_j}{\|Q_j x^r\|}~~\mbox{for all}~r \in R'.
\]
Taking the limit as $r \rt \infty$ shows
\[
s \ge \frac{w_j}{\|Q_j \bar{x}\|} \ge \min_i \frac{w_i}{\|Q_i \bar{x}\|},
\]
as required.  Thus the function $\Phi$ is indeed lower semicontinuous.

Denote the right hand side of the second claimed expression for $\Phi$
by $\Psi((y_i))$:  we next want to prove that this infimum is attained
whenever $\Psi((y_i))$ is finite.
Notice that the infimum is unchanged if we
add the condition $\|y^*\| = 1$, using positive homogeneity.  Now
suppose that the infimum is finite, so there exist feasible
vectors $\bar{u}_i^*$ and $\bar{y}^*$.  If we define 
$\beta = \max_i \|\bar{u}_i^*\|$, then we can rewrite the infimum as
\begin{eqnarray*}
\lefteqn{
\inf_{u_i^* \in U_i^*,\:y^* \in Y^*}
\Big\{
\max_i \|u_i^*\| : \sum_i \ip{y^*}{y_i} Q_i^* u_i^* \in F^*(-y^*),~
\|y^*\| = 1,
} \\
 & & \mbox{} \hspace{8cm} \ip{y^*}{y_i} \ge 0,~\|u_i^*\| \le \beta~\forall i
\Big\}.
\end{eqnarray*}
This is the infumum of a continuous function over a nonempty compact set,
so is attained.

It remains to prove that the two functions $\Phi$ and $\Psi$ are identical.
Consider any real $\psi > 0$.
Using the attainment property we have just proved for $\Psi$,
the statement $\Psi((y_i)) \le \psi$ is
equivalent to the solvability of the system
\begin{eqnarray*}
\sum_i \ip{y^*}{y_i} Q_i^* u_i^* \in F^*(-y^*), & & 0 \ne y^* \in Y^* \\
\ip{y^*}{y_i} \ge 0, & & u_i^* \in \psi B_{U_i^*}~~\mbox{for each}~i,
\end{eqnarray*}
or equivalently, to the solvability of the system
\begin{eqnarray*}
\sum_i \ip{y^*}{\psi^{-1} y_i} Q_i^* u_i^* \in F^*(-y^*), & & 0 \ne y^* \in Y^* \\
\ip{y^*}{\psi^{-1} y_i} \ge 0, & & u_i^* \in B_{U_i^*}~~\mbox{for each}~i.
\end{eqnarray*}
Using the theorem of the alternative (\ref{4.1}), this is equivalent to the
unsolvability of the system
\[
\sum_i w_i \psi^{-1} y_i \in F(x),~~
\|Q_i x\| < w_i \in \R~\mbox{for each}~i,~~x \in X,
\]
or equivalently (since $F$ is positively homogeneous), to the unsolvability
of the system
\[
\sum_i w_i y_i \in F(x),~~
\psi < \frac{w_i}{\|Q_i x\|},~~0 < w_i \in \R~\mbox{for each}~i,~~x \in X.
\]
But this in turn is equivalent to the statement $\Phi((y_i)) \le \psi$.
To summarize, we have shown, for all real $\psi > 0$,
\[
\Psi \Big( (y_i) \Big) \le \psi
~~\Leftrightarrow~~
\Phi \Big( (y_i) \Big) \le \psi.
\]
The result now follows.
\finpf

\section{The main result}
We now have all the tools we need to derive our main result.

\begin{thm}[distance to nonsurjectivity]
For any finite-dimen\-sional normed spaces $X,Y,U_i,V_i$, 
closed sublinear surjective set-valued 
mapping $F:X \tto Y$, and linear mappings $P_i:V_i \rt Y$ and $Q_i:X \rt U_i$
(for $i=1,2,\ldots,k$), the following four quantities are equal:
\[
\inf_{T_i \in L(U_i,V_i)}
\Big\{
\max_i \|T_i\| : F + \sum_i P_i T_i Q_i~\mbox{\rm nonsurjective}
\Big\};
\]
\[
\inf_{\mbox{\scriptsize\rm rank-one}~T_i \in L(U_i,V_i)}
\Big\{
\max_i \|T_i\| : F + \sum_i P_i T_i Q_i~\mbox{\rm nonsurjective}
\Big\};
\]
\[
\inf_{u_i^* \in B_{U_i^*},\:z_i \ge 0,\:0 \ne y^* \in Y^*}
\Big\{
\max_i \frac{z_i}{\|P_i^* y_i\|} : 
\sum_i z_i Q_i^* u_i^* \in F^*(y^*)
\Big\};
\]
\[
\inf_{v_i \in B_{V_i}}~
\sup_{x \in X,\:w_i > 0}
\Big\{
\min_i \frac{w_i}{\|Q_i x\|} : \sum_i w_i P_i v_i \in F(x)
\Big\}.
\]
Furthermore, if these quantities are finite, each infimum above is
attained.
\end{thm}

\pf
The equality of the first three expressions follows immediately from
Theorem \ref{3.5n} (rank reduction and surjectivity).  The last expression
also follows from the same result, after applying the duality
theorem (\ref{4.7}).
\finpf


\begin{thebibliography}{10}

\bibitem{Bor86}
J.M. Borwein.
\newblock Norm duality for convex processes and applications.
\newblock {\em Journal of Optimization Theory and Applications}, 48:53--64,
  1986.

\bibitem{Don03}
A.L. Dontchev, A.S. Lewis, and R.T. Rockafellar.
\newblock The radius of metric regularity.
\newblock {\em Transactions of the American Mathematical Society},
  355:493--517, 2003.

\bibitem{Doy82}
J.~Doyle.
\newblock Analysis of feedback systems with structured uncertainty.
\newblock {\em IEEE Preceedings}, 129:242--250, 1982.

\bibitem{Lew99Ill}
A.S. Lewis.
\newblock Ill-conditioned convex processes and linear inequalities.
\newblock {\em Mathematics of Operations Research}, 24:829--834, 1999.

\bibitem{Pen98}
J.~{Pe\~{n}a}.
\newblock {\em Condition numbers for linear programming}.
\newblock PhD thesis, Cornell University, 1998.

\bibitem{Pen00}
J.~{Pe\~{n}a}.
\newblock Understanding the geometry of infeasible perturbations of a conic
  linear system.
\newblock {\em SIAM Journal on Optimization}, 10:534--550, 2000.

\bibitem{Pen03}
J.~{Pe\~{n}a}.
\newblock A characterization of the distance to infeasibility under structured
  perturbations.
\newblock {\em Linear Algebra and its Applications}, 2003.
\newblock To appear.

\bibitem{Ren95a}
J.~Renegar.
\newblock Incorporating condition measures into the complexity theory of linear
  programming.
\newblock {\em SIAM Journal on Optimization}, 5:506--524, 1995.

\bibitem{Ren95}
J.~Renegar.
\newblock Linear programming, complexity theory and elementary functional
  analysis.
\newblock {\em Mathematical Programming}, 70:279--351, 1995.

\bibitem{Rob76}
S.M. Robinson.
\newblock Regularity and stability for convex multivalued functions.
\newblock {\em Mathematics of Operations Research}, 1:130--143, 1976.

\bibitem{Roc70}
R.T. Rockafellar.
\newblock {\em Convex Analysis}.
\newblock Princeton University Press, Princeton, N.J., 1970.

\bibitem{Roc98}
R.T. Rockafellar and R.J.-B. Wets.
\newblock {\em Variational Analysis}.
\newblock Springer, Berlin, 1998.

\end{thebibliography}

\end{document}